# A rigorous lower confidence bound for the expectation of a positive random variable

Yoram Gat

*Intel Labs*
*SC12-303*
*2200 Mission College Blvd.*
*Santa Clara, CA 95054-1537*
*e-mail:* yoram.gat@intel.com

**Abstract:** Given an IID sample from a positive distribution, we provide a method for constructing rigorous finite sample lower confidence bounds for the expectation of the distribution. The method is based on constructing rigorous confidence regions for the cdf of the distribution. We provide some analysis of the asymptotical behavior of the rigorous LCBs. We apply the method to obtain an LCB for a particular, controversial, empirical data set, where the validity of standard methods has been called into question.

## 1. Introduction

In general, given an IID sample from an unknown distribution, no rigorous confidence bounds can be provided for the expectation of the distribution. The possibility of the existence of a very low probability tail with some extreme value which has a large impact on the expectation can never be ruled out, or even made improbable by any finite sample. For example, with a sample of size $n$, the existence of an atom with probability $n^{-2}$ and an arbitrary value, with magnitude large enough so as to impact the expected value greatly, is not unlikely. Indeed, it may be that an expectation does not even exist, or has infinite magnitude.

The normal practice for estimating expectations is to ignore the possibility of the existence of low probability, extreme value tails, apply estimators with known asymptotical properties, and, in effect, assume that those properties are valid for the given sample size. Alternatively, rigorous confidence bounds, such as the Chebychev or Chernoff bounds, can be derived when certain moments are assumed to be bounded or when the distribution itself is assumed to lie within a bounded interval.

Here, we use a weaker assumption, namely, that the distribution is over positive numbers only, but aim at deriving only lower confidence bounds rather than a confidence interval. The same argument as above shows that without additional assumptions, no upper confidence bounds for the expectation can be established.

In addition to providing a guaranteed finite sample confidence level for any underlying distribution over the positive reals, the confidence bounds proposed have the pleasing property that they are monotonic in the order statistics of the sample. This eliminates the paradoxical phenomenon, which can occur with a normal-theory LCB, in which a positive outlier in the sample lowers the nominal LCB for the expectation.

## 2. Setup and theorem

Let $X, X_1, \ldots, X_n$ be IID variables from an unknown distribution over the real numbers, with $\mathbf{P}(X \geq 0) = 1$. Let $X_{(1)}, \ldots, X_{(n)}$ be the order statistics of $X_1, \ldots, X_n$. We wish to derive lower confidence bounds for the expectation of $X$: $B = B(X_1, \ldots, X_n)$ such that $\mathbf{P}(B > \mathbf{E}X) \leq \alpha$. We rely on the following theorem. The theorem establishes a LCB for $\mathbf{E}X$ as a consequence of simultaneous UCBs for the cdf of $X$ at the sample points.

**Theorem 1.** *Let $\mathbf{U} = \big(U_{(1)}, \ldots, U_{(n)}\big)$ be the vector of the order statistics of $n$ independent samples from a $\mathrm{U}[0,1]$ distribution. For any vector $u = (u_1, \ldots, u_n) \in [0,1]^n$, define*
$$p_u = \mathbf{P}(U_{(1)} \leq u_1, \ldots, U_{(n)} \leq u_n).$$

*Define*
$$B_u = \sum_{k=1}^n (1 - u_k)\big(X_{(k)} - X_{(k-1)}\big) \tag{1}$$
$$= \sum_{k=1}^n (u_{k+1} - u_k) X_{(k)}, \tag{2}$$

*where $u_{n+1} = 1$ and $X_{(0)} = 0$. Then $B_u$ is a level-$p_u$ LCB for the expectation of $X$.*

*Proof:* The second equality follows from rearranging the terms of the sum. We prove the first equality: Let $F^-(x)$ be the left-continuous cdf of the random variable $X$, i.e.,
$$F^-(x) = \mathbf{P}(X < x).$$
Then for any monotonically non-decreasing vector $x = (x_1, \ldots, x_n) \in \mathbb{R}^{+n}$,
$$\sum_{k=1}^n \big(1 - F^-(x_k)\big)(x_k - x_{k-1}) \leq \mathbf{E}X,$$
where $x_0 = 0$.

The random variables $F^-(X_1), \ldots, F^-(X_n)$ are independent, identically distributed variables such that $P(F^-(X_i) \leq x) \geq x$. Therefore, there exists a set of IID $U[0,1]$ random variables $U_1, \ldots, U_n$, such that $U_i \geq F(X_i)$, $i = 1, \ldots, n$.

It follows that for any vector $u = (u_1, \ldots, u_n) \in [0,1]^n$ the probability of the event
$$\mathbf{P}(F(X_{(1)}) \leq u_1, \ldots, F(X_{(n)}) \leq u_n)$$
is at least $p_u$. On this event,
$$B_u \leq \sum_{k=1}^n \left(1 - F^-(X_{(k)})\right)\left(X_{(k)} - X_{(k-1)}\right) \leq \mathbf{E}X. \quad \square$$

## 3. Tunable families of bound parameter vectors

Theorem 1 implies that each vector $u \in [0,1]^n$ defines a level-$p_u$ LCB, $B_u$, for $\mathbf{E}X$. Different choices of the parameter vector $u$ result in different LCBs. It is convenient to construct families of parameter vectors in such a way that from each family a vector can be chosen to match a desired level of confidence: Let $\Lambda$ be a closed subset of $\mathbb{R}$. Define a tunable family of parameter vectors, $\mathcal{U} = \{u^\lambda \in [0,1]^n : \lambda \in \Lambda\}$, to be a set of vectors which is parameterized continuously by $\lambda$ and increasing monotonically in each coordinate to 1. That is, for all $i = 1, \ldots, n$, the following hold:

- If $\lambda_1 \leq \lambda_2$ then $u_i^{\lambda_1} \leq u_i^{\lambda_2}$.
- If $\lim_k \lambda_k = \lambda$ then $\lim_k u_i^{\lambda_k} = u_i^\lambda$.
- $\sup_{\lambda \in \Lambda} u_i^\lambda = 1$.

Then for each $\alpha$, $0 < \alpha < 1$, there exists a unique $\lambda(\alpha)$ such that
$$\lambda(\alpha) = \min\left\{\lambda \in \Lambda : p_{u^\lambda} \geq 1 - \alpha\right\}.$$

Given a desired confidence level $1 - \alpha$, or a set of confidence levels $1 - \alpha_1, 1 - \alpha_2, \ldots, 1 - \alpha_k$, the corresponding $\lambda(\alpha)$, or $\lambda(\alpha_1), \lambda(\alpha_2), \ldots, \lambda(\alpha_k)$, can be determined numerically, to any desired precision, using simulation.

**Examples** An infinite variety of tunable parameter vector families exist. A very simple family is defined by adding an offset to a the vector $(\frac{1}{n+1}, \ldots, \frac{n}{n+1})$:
$$\mathcal{U}_{\text{OFF}} = \left\{\left(\min(1, \frac{1}{n+1} + \lambda), \ldots, \min(1, \frac{n}{n+1} + \lambda)\right) : \lambda \in [0,1]\right\}.$$

When using this family, the LCB assumes the form
$$\text{LCB}_{\text{OFF}} = \frac{1}{n+1} \sum_{i=1}^{n_\lambda - 1} X_{(i)} + \left(1 - \lambda - \frac{n_\lambda - 1}{n+1}\right) X_{(n_\lambda)},$$
where $n_\lambda = \lceil (n+1)(1-\lambda) \rceil$.

Another family results from using confidence bounds for the beta distribution. Let $B(a, b, p)$ be a level-$p$ UCB for the beta distribution with parameters $a$ and $b$. The construction of this family is intuitively motivated by the fact that for continuous IID random variables $X_1, \ldots, X_n$ the marginal distribution of $F(X_{(i)})$ is a beta distribution with parameters $i$ and $n-i+1$ for all $i = 1, \ldots, n$. Define:

$$\mathcal{U}_{\text{BETA}} = \left\{ \Big(B(1, n, \lambda), B(2, n-1, \lambda), \ldots, B(n, 1, \lambda)\Big) : \lambda \in [0, 1] \right\}.$$

## 4. Some asymptotical analysis

Analysis of the asymptotical behavior of the rigorous LCB is facilitated by the Donsker property of the empirical process (see, for example, [1], chapter 2.1). The Donsker property implies that for continuous distributions the centered and scaled empirical process converges in distribution to the standard Brownian bridge. The centered and scaled empirical process, $(H_n(t), 0 \leq t \leq 1)$, is defined as:

$$H_n(t) = \frac{1}{\sqrt{n}} \left( \sum_{i=1}^{n} \mathbf{1}(X_i \leq t) - t \right).$$

This property can be used directly to calculate the asymptotical behavior of the rigorous bound obtained when using the offset family, $\text{LCB}_{\text{OFF}}$, as we do below. Asymptotical analysis of the rigorous LCB for other families such as $\text{LCB}_{\text{BETA}}$ would be more complex.

The distribution of $M$, the supremum of the Brownian bridge, is

$$\mathbf{P}(M > x) = e^{-2x^2},$$

putting the $1 - \alpha$ quantile of the distribution, $q_\alpha$ at $\sqrt{\frac{1}{2} \log \frac{1}{\alpha}}$. When using the offset family with a sample of size $n$, the member selected for a $1 - \alpha$ LCB will be approximately $(\frac{1}{n+1} + \frac{q_\alpha}{\sqrt{n}}, \ldots, \frac{n}{n+1} + \frac{q_\alpha}{\sqrt{n}})$. If the distribution of $X$ has a finite mean $\text{LCB}_{\text{OFF}}$ is equal to

$$\frac{1}{n+1} \sum_{i=1}^{n - \lceil q_\alpha \sqrt{n} \rceil} X_{(i)} + o\left(\frac{1}{\sqrt{n}}\right),$$

its expected value is

$$\mathbf{E}\,\text{LCB}_{\text{OFF}} = \mathbf{E} X \mathbf{1}\left\{ X \leq F_X^{-1}\left(1 - \frac{q_\alpha}{\sqrt{n}}\right) \right\} + o\left(\frac{1}{\sqrt{n}}\right),$$

and so,

$$\mathbf{E} X - \mathbf{E}\,\text{LCB}_{\text{OFF}} = \mathbf{E} X \mathbf{1}\left\{ X > F_X^{-1}\left(1 - \frac{q_\alpha}{\sqrt{n}}\right) \right\} + o\left(\frac{1}{\sqrt{n}}\right).$$

The first term on the right is the integral of the tail of the distribution. If $\mathbf{E} X$ is finite, this term approaches zero as $n$ increases, guaranteeing that the LCB is consistent. However, the convergence may be arbitrarily slow unless some additional assumptions regarding the distribution of $X$ are made. The convergence is $O(n^{-\frac{1}{2}})$ if and only if $X$ is bounded almost surely. Using the Hölder inequality it can be shown that if $\mathbf{E}X^{r+\epsilon} < \infty$ for some positive $r$ and $\epsilon$ then the convergence is $o(n^{-\frac{1}{2}+\frac{1}{2r}})$.

Of course, for the normal theory LCB, $\text{LCB}_{\text{Normal}}$,

$$\sqrt{n}\,(\mathbf{E}X - \mathbf{E}\,\text{LCB}_{\text{Normal}}) = \Phi(1-\alpha/2)\sqrt{\mathbf{Var}X},$$

guaranteeing $O(n^{-\frac{1}{2}})$ if $\mathbf{Var}X$ is finite, but no convergence if $\mathbf{Var}X$ is infinite.

Thus, asymptotically, the normal theory LCB converges faster than the rigorous LCB whenever $\mathbf{Var}X$ exists (unless $X$ is bounded, in which case both LCBs converge as $O(n^{-\frac{1}{2}})$), but the rigorous LCB guarantees convergence to $\mathbf{E}X$ when $\mathbf{Var}X$ is infinite, i.e., in situations in which the normal-theory LCB diverges in expectation.

## 5. Application to the Lancet study of mortality in Iraq

The results above provide a method for generating theoretical rigorous lower confidence bounds for the expectation of a positive random variable. These bounds can be applied in situations where conventional methods for producing confidence bounds are challenged based on the fact that the validity of those methods relies on asymptotical analysis which may not hold for samples of a given size from particular a distribution.

One such case is the politically sensitive estimate of mortality in Iraq following the U.S. led invasion. In 2006 a group of researchers from the Johns Hopkins Bloomberg School of Public Health carried out a survey among households in Iraq aimed at estimating mortality [2]. They provided a point estimate of about 601,000 violent deaths in Iraq for the period March 2003 to July 2006 and a 95% confidence interval of 426,000-794,000. Due to the potential political implication of the findings, the study received intense scrutiny. Most of the attention was focused at the various potential biases introduced into the data collection process by a methodology constrained by the conditions in Iraq (see a summary of such points of criticism in [3]). In addition, however, and despite the fact that the estimation procedure used was apparently identical to that used in similar studies, some criticism was made of the estimation procedure itself. Doubts were voiced as to whether the normal theory 95% confidence interval did indeed have its nominal probability of coverage and it was suggested that the a true 97.5% LCB would be drastically lower than the left point of the interval [4].

We follow here the treatment of Mark van der Lann [4]. He uses a somewhat stylized setup in which the death counts in the 49 clusters in the sample collected by Burnham et al. are assumed to be IID samples from an unknown distribution of violent deaths in household clusters in Iraq[1]. Each cluster contains 40 households, so under van der Laan's setup the unknown mean of the distribution is 40 times the mean number of violent deaths per household in Iraq.

---

[1] In reality, the sample was stratified geographically by governorates.

Van der Lann provides the death counts in the 47 clusters as follows:

| 17 | 15 | 0 | 0 | 3 | 4 | 0 | 0 |
|----|----|---|---|---|----|----|---|
| 5 | 1 | 0 | 2 | 2 | 6 | 0 | 9 |
| 3 | 0 | 1 | 5 | 7 | 12 | 22 | 0 |
| 0 | 0 | 1 | 25 | 2 | 24 | 35 | 9 |
| 6 | 5 | 4 | 4 | 6 | 1 | 3 | 1 |
| 3 | 5 | 2 | 0 | 25 | 9 | 18 | |

The sample mean is 6.4 and the sample standard deviation is 8.3, giving a classical normal-theory 97.5% LCB for the expectation of the distribution of 4.0, i.e., 63.0% of the sample mean. Employing the method above, we obtain a rigorous LCB for the expectation of 2.3 (36.5% of the sample mean) using the offset family and 2.8 (43.8% of the sample mean) when using the beta family.[2] We therefore note that while the rigorous LCBs constructed here are significantly lower than the nominal normal-theory LCB, they are not dramatically different (reducing the bound by about one third). This suggests that using such a technique can be useful when dealing with certain situations in which the validity of traditional methods may be called into doubt.

Figure 1 demonstrates the construction of the LCBs graphically. It shows the empirical cdf together with lines signifying the boundaries of the confidence regions established for the cdf using the offset family (dotted line) and the beta family (dashed line). The LCBs for the expectation are the areas to the left of and above those two curves.

## 6. Further research

One point associated with the method presented that may merit further research regards the choice of tunable family. Are some families better - i.e., yield tighter bounds - than others, across all possible distributions? Can families be chosen so as to match various properties of the distribution?

Another avenue of research would be to produce extensions of the method in order to make it applicable to a wider variety of situations. One desirable extension would be to cover cases where the sample is stratified, while another would be to cases where random censoring occurs.

---

[2]Employing a simple linear relationship between the expectation of deaths in a cluster and the total number of deaths in the population, these bounds would correspond to LCBs for the total number of violent Iraqi deaths of 219,000 (offset family) and 263,000 (beta family) in the period covered by the survey.

**Iraq mortality empirical cdf and confidence bounds**

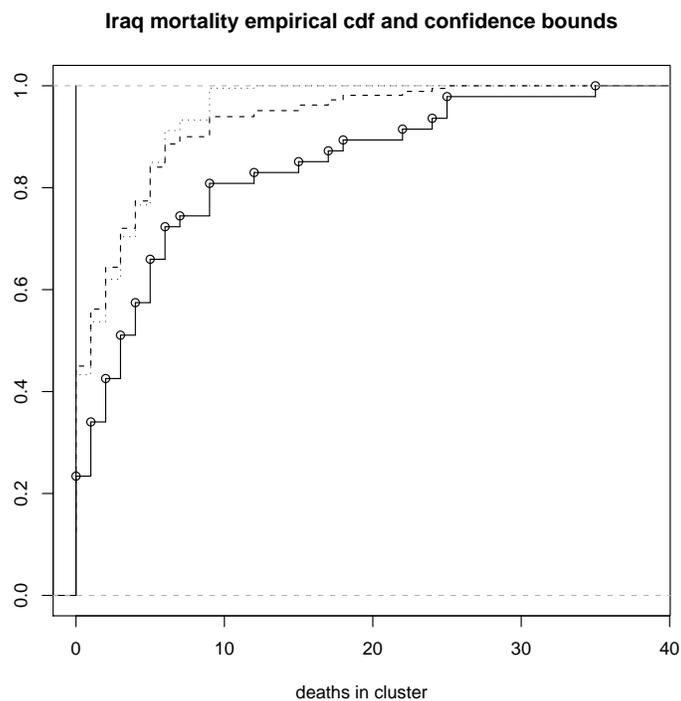

FIG 1. *The Burnham et al. cluster death counts empirical cdf, and derived confidence regions.*